# A New Efficient Stochastic Energy Management Technique for Interconnected AC Microgrids

Morteza Dabbaghjamanesh, *Student Member, IEEE,* Shahab Mehraeen, *Member, IEEE*
Abdollah Kavousi-Fard, *Member, IEEE,* Farzad Ferdowsi, *Member, IEEE*

***Abstract*— Cooperating interconnected microgrids with the Distribution System Operation (DSO) can lead to an improvement in terms of operation and reliability. This paper investigates the optimal operation and scheduling of interconnected microgrids highly penetrated by renewable energy resources (DERs). Moreover, an efficient stochastic framework based on the Unscented Transform (UT) method is proposed to model uncertainties associated with the hourly market price, hourly load demand and DERs output power. Prior to the energy management, a newly developed linearization technique is employed to linearize nodal equations extracted from the AC power flow. The proposed stochastic problem is formulated as a single-objective optimization problem minimizing the interconnected AC MGs cost function. In order to validate the proposed technique, a modified IEEE 69-bus network is studied as the test case.**

## NOMENCLATURE

***Indices***

| | |
|---|---|
| $d$ | Index for loads |
| $i$ | Index of DERs |
| $j$ | Index for adjacent MGs |
| $k$ | Index of MG |
| $M$ | Index of main grid |
| $n,m$ | Indices of nodes in feeder |
| $t$ | Index of time |

***Sets***

| | |
|---|---|
| $B$ | Set of MGs |
| $BMG$ | Set of adjacent MGs/utility |
| $D$ | Set of adjustable loads |
| $G$ | Set of dispatchable units |
| $\eta$ | Set of remaining nodes |

***Parameters***

| | |
|---|---|
| $E$ | Adjustable load required energy |
| $C,C^M$ | Cost, Cost of adjacent MG/main grid ($/kW) |
| $DT$ | Minimum Down time |
| $RD$ | Ramp down rate |
| $RU$ | Ramp up rate |
| $S_p$ | Constant power loads |
| $UT$ | Minimum Up time |
| $Y_{nm}$ | Branch admittance through buses *n, m* |

***Variables***

| | |
|---|---|
| $I$ | Current |
| $I_{slack}$ | Slack node current |
| $P$ | Active power |
| $P^M$ | Main grid power |
| $Q$ | Reactive power |
| $SU,SD$ | Start up, shot down costs |
| $SP_{nm}$ | Thermal limit of line *nm* |
| $T^{on}$ | Number of successive ON hours |
| $T^{off}$ | Number of successive OFF hours |
| $U$ | DER status |
| $V_n$ | Voltage at bus n |
| $W^0$ | Weight of the mean value |
| $\mu$ | Mean value |
| $\lambda$ | Covariance |
| $\delta_n$ | Voltage angle of the *nth* bus |
| $\theta_{nm}$ | Angle of branch through buses *n, m* |

M.dabbaghjamanesh, S. Mehran and F. Ferdowsi are with the Department of Electrical and Computer Engineering, Louisiana State University, Baton Rouge, LA, 70803, USA. (E-mail: mdabba1@lsu.edu, smehraeen@lsu.edu, fferdowsi1@lsu.edu).
A. Kavousi-Fard is with the Department of Electrical and Computer Engineering, University of Michigan, MI, 48126, USA. (E-mail: akavousi@umich.edu).

## I. INTRODUCTION

MICROGRIDS (MGs) have recently emerged into the traditional AC networks due to their technical and economic benefits in presence of high penetration of renewable energy resource (RESs) either in grid-connected or islanded mode [1]. Significant benefits of MGs incorporate both technical and economic aspects such as higher reliability and resiliency, self-healing capabilities and lower cost in operational and design stages. Higher power quality, lower losses, less planning and operation costs, higher reliability and avoiding market monopoly are some of the other advantages of MGs.

In traditional power grids, the authority is given to the Distribution System Operation (DSO) to plan and manage the distribution feeders. However, in modern power systems, the DSO and MG may have different utility operators. Therefore, power/energy management in different microgrids may be based on different rules and policies. Since the whole entire grid is made up of interconnected interconnected MGs, any change in any of subsystems affects the network operation from the system level perspective due to the high level of interconnectivity within the system. According to the IEEE-1547.4 standard [2], decoupling the network into interconnected sub-MG systems can improve the system operation and reliability significantly. However, there might be some limitations on communications among different MGs due to the privacy and security issues. Because of these constraints, a system-level coordination method is required to accomplish an efficient energy management for interconnected MGs while different perturbations occur in the system such as load switching or changes in generation units.

In different sets of literature, planning, operation, energy management, and service restoration are studied for MGs [3-8]. In [3], a robust optimization technique is proposed for MG planning by taking the load and DG generation uncertainties into account. An efficient energy management for multi-period islanding of MGs is studied in [4]. Authors have investigated the energy management problem utilizing bender decomposition technique in grid-connected and islanded operation modes. A two-layer dispatch framework is presented in [5-8] to solve the economic dispatch problem considering the spinning reserve. The market-based MG operation scheduling is studied in [9] when the MG is in correlation with the distribution market operator. Authors have used the Mixed-Integer Linear Programming (MILP) technique to solve the problem. In [10] and [11], the MG operation is studied when the network is highly penetrated by RESs along with storage units utilized to overcome uncertainties in generation units.

Although the energy management problem in MGs is widely discussed in different research studies, the impact of interconnectivity within a system made up of interconnected MGs is not well discussed yet. In [12] and [13], the generation scheduling along with the mismatch control of multi-agents interconnected MGs is studied. In [14] the energy management of interconnected MGs in the grid-connected mode is investigated when the demand is considered as an unknown variable for some adjacent MGs. A bi-level programming-based energy management of interconnected MGs is proposed in [15] where the upper and lower levels are central generation units and DERs, respectively. Optimal power dispatch of interconnected MGs is addressed in [16] considering uncertainties in generation and demand. A robust optimization method is studied in [17] for optimal scheduling of renewable-based distributed MGs. A decentralized energy management technique for interconnected MGs based on Alternating Direction Method of Multipliers (ADMM) method is presented in [18]. Authors have modeled uncertainties associated with RESs and load demand utilizing Monto Carlo Simulation (MCS) followed by the scenario reduction technique.

While previous research studies have addressed interconnected MGs operation challenges from different aspects, the stochastic-based energy management of interconnected MGs can be improved which is only discussed in a few sets of literature. This paper investigates the optimal energy management of interconnected MGs penetrated by intermittent distributed energy resources (DERs). Due to the stochastic nature of the proposed problem, Unscented Transform (UT) technique is utilized to model parameter uncertainties associated with hourly load demand, hourly market price, and RESs' output generation. The feasibility and satisfying performance of the proposed stochastic framework are examined on an IEEE test syste. To summarize, the main contributions of this paper can be summarized as follows:

- Providing a sufficient mathematical model for optimal energy management in interconnected AC MGs. The proposed formulation enables centralized processing with the minimal exchange of information among MGs.
- Applying a newly developed linearization approach for AC power flow equations in the power distribution network level.
- Developing a new UT-based stochastic framework for interconnected MGs operating under uncertainties such as hourly load demand, hourly market price, and RESs' output power.

## II. Interconnected Microgrid Scheduling Formulation

### *A. Objective Function and Problem Constraints*

The cost function includes the costs associated with the power purchased from the adjacent MG/main grid and the cost imposed by the local utility (generation units within each MG):

$$\min \sum_{k\in B}\sum_{t}\sum_{i\in G}[C_{ikt}P_{ikt}+SU_{ikt}+SD_{ikt}]+\sum_{j\in BMG}C^{M}_{jkt}P^{M}_{jkt} \qquad \forall t \qquad (1)$$

The optimization problem includes the essential constraints described in (2) -(14):

*Maximum power exchanged between the MGs and the main grid:* Equation (2) shows the constraint on the amount of power exchanged between each two MGs or a MG and the main grid:

$$-P^{M,\min}_{jk} \le P^{M}_{jkt} \le P^{M,\max}_{jk}. \qquad \forall t, \forall k\in B, \forall t \qquad (2)$$

*DG Generation Capacity:* Each DG has a limitation in terms of power generation described in (3) and (4):

$$P_{ikt,\min}U_{ikt} \le P_{ikt} \le P_{ikt,\max}U_{ikt} \qquad \forall i\in G, \forall k\in B, \forall t \qquad (3)$$

*Ramp up/down constraints:* Equations (4) and (5) represent the ramp up/down rate constraints on the generation units:

$$P_{ikt}-P_{ik(t-1)} \le RU_{ik} \qquad \forall i\in G, \forall k\in B, \forall t \qquad (4)$$

$$P_{ik(t-1)}-P_{ikt} \le RD_{ik} \qquad \forall i\in G, \forall k\in B, \forall t \qquad (5)$$

*Minimum up and down times constraints:* The dispatchable generation units are subject to a minimum and maximum up and down times limits given in (6) and (7):

$$T^{on}_{ikt} \ge UT_{ik}(U_{ikt}-U_{ik(t-1)}) \qquad \forall i\in G, \forall k\in B, \forall t \qquad (6)$$

$$T^{off}_{ikt} \le DT_{ik}(U_{ik(t-1)}-U_{ikt}) \qquad \forall i\in G, \forall k\in B, \forall t \qquad (7)$$

*Energy storage constraints:* The battery state of charge (SOC) must remain in a permissible limit in all time intervals. Equation (8) represents the battery's SOC constraint:

$$P_{uk,\min} \le P_{ukt} \le P_{uk,\max}. \qquad \forall k\in B, \forall t \qquad (8)$$

*Adjustable loads constraint:* Adjustable loads must be within their limits as follows:

$$D_{dkt,\min} \le D_{dkt} \le D_{dkt,\max} \qquad \forall d\in D, \forall k\in B, \forall t \qquad (9)$$

However, at the end of the horizon the required energy of each adjustable load must be satisfied as follows:

$$\sum D_{dkt} = E_{dk} \qquad \forall d\in D, \forall k\in B \qquad (10)$$

### B. Network Constraints

*AC power flow constraints:* The nodal equations described in (11) represent the power balance for the Ac grid:

$$\begin{cases} \sum_{k\in B} P_{nkt} = \sum_{m} |V_{nkt}||V_{mkt}||Y_{nmk}| \times \cos\left(\theta_{nmk} + \delta_{nkt} - \delta_{mkt}\right) \\ \sum_{k\in B} Q_{nkt} = \sum_{m} |V_{nkt}||V_{mkt}||Y_{nmk}| \times \sin\left(\theta_{nmk} + \delta_{nkt} - \delta_{mkt}\right) \end{cases} \quad (11)$$

*Bus voltage limit:* The steady-state voltage of each bus must be within the limits as described in (12) and (13) in terms of the voltage stability:

$$V_{nk\,\min} \le V_{nkt} \le V_{nk\,\max} \quad (12)$$

$$-\pi \le \theta_{nkt} \le \pi \quad (13)$$

*Maximum power flow in feeders:* The power flow in any feeder is not expected to violate its maximum permissible value described as an inequality constraint in (14).

$$SP_{nmtk} \le SP_{nmk\,\max} \quad (14)$$

## III. Load Flow Analysis for Power Distribution System

Linearized power flow analysis in the form of an approximated DC power flow is commonly used for AC networks. However, due to the high *R/X* ratio in the distribution network, this can lead to inaccurate results. In this paper a newly developed linearization approximation technique is utilized in which the radial topology is not necessarily required [19].

Based on the network's nodal equations, the voltages and currents are related to the grid's admittance matrix as shown in (15) [19].

$$\begin{pmatrix} I_{slack} \\ I_{\eta} \end{pmatrix} = \begin{pmatrix} Y_{slack,slack} & Y_{slack,\eta} \\ Y_{\eta,slack} & Y_{\eta,\eta} \end{pmatrix} \cdot \begin{pmatrix} V_{slack} \\ V_{\eta} \end{pmatrix} \quad (15)$$

In addition, based on the ZIP load model, the current at each node can be represented as a combination of three terms as described in (16) [19]:

$$I_{\eta} = \frac{S^*_{P\eta}}{V^*_{\eta}} + h S^*_{I\eta} + h^2 S^*_{Z\eta} V_{\eta} \quad (16)$$

where $h=1/V_{norm}$ (per unit). It should be noted that the ZIP configuration is a linear model except at the constant power load part. However, this term is approximated to achieve a linear power flow [19].

By approximating the voltage as $V=(1-\Delta V)$ and developing the Taylor series, a linear form of the voltage can be obtained as follows if high order terms are neglected [19]:

$$\frac{1}{V} = \frac{1}{1-\Delta V} \approx 1 + \Delta V = 2 - V \quad (17)$$

The linear representation of the ZIP model described in (16) is formed in (18) is obtained as follows [19]:

$$I_{\eta} = h S^*_{P\eta}.(2 - hV^*_{\eta}) + h S^*_{I\eta} + h^2 S^*_{Z\eta} V_{\eta} \quad (18)$$

The load model is reformed in (19) by some manipulations [19]:

$$A_1 + A_2 V^*_{\eta} + A_3 V_{\eta} = 0 \quad (19)$$

where,

$$A_1 = Y_{\eta,slack} V_{slack} - 2h S^*_{P\eta} - h S^*_{I\eta}$$

$$A_2 = h^2 dig(S^*_{P\eta})$$

$$A_3 = Y_{\eta\eta} - h^2 dig(S^*_{Z\eta})$$

## IV. A Stochastic Optimization Framework Based on the Unscented Transform

As that is already mentioned, the proposed optimization problem contains some uncertainties in parameters such as hourly market price, hourly load demand, and output power from RESs. The Unscented Transform (UT) technique is one of the well-known estimation methods bringing significant advantages into the analysis such as easy coding, low computational burden and high ability for modeling the nonlinear uncertainties [20]. This method is explained in the rest:

A nonlinear function $y$=f($x$), with $\alpha$ random inputs and $h$ random outputs with the mean value $\mu$ and covariance $\lambda$ is considered. According to the UT method, the mean value and the covariance of the output ($\mu_y$ and $\lambda_y$ respectively) have to be calculated as follows [20]:

Step 1: Calculate $2\alpha+1$ samples from the input data described in (20) -(22).

$$x_0 = \mu \quad (20)$$

$$x_{\omega} = \mu + \left(\sqrt{\frac{\alpha}{1-W^0}\lambda}\right)_{\omega} \quad , \omega = 1,2,...,\alpha \quad (21)$$

$$x_{\omega} = \mu - \left(\sqrt{\frac{\alpha}{1-W^0}\lambda}\right)_{\omega} \quad , \omega = 1,2,...,\alpha \quad (22)$$

Step 2: calculate the weighting factor of each sample point as presented in (23) -(26).

$$W^0 = W^0 \quad (23)$$

$$W_{\omega} = \frac{1-W^0}{2\alpha}; \quad \omega = q+1,...,\alpha \quad (24)$$

$$W_{\omega+n} = \frac{1-W^0}{2\alpha}; \quad \omega+\alpha = \alpha+1,...,2\alpha \quad (25)$$

$$\sum W_{\omega} = 1 \quad (26)$$

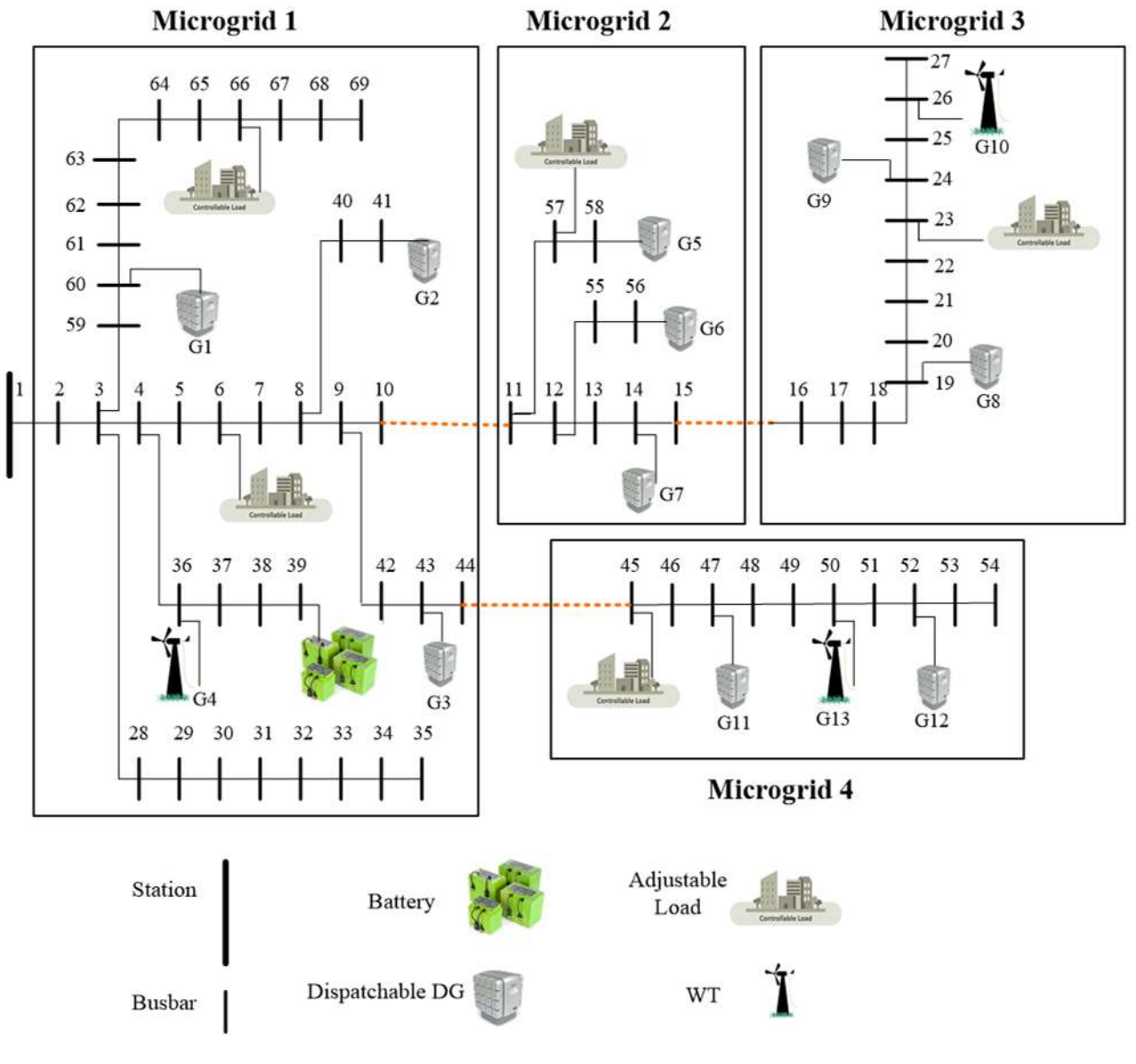

Fig. 2. One-line diagram of the test interconnected microgrid

Step 3: Using the nonlinear function $y$=f($x$) to find the output samples as follows:

$$y_{\omega} = f(X_{\omega}) \quad (27)$$

Step 4: Calculate the covariance and the mean value ( $\mu_y$ and $\lambda_y$ ) for the output variable £ as follows:

$$\mu_y = \sum_{\omega} W_{\omega} £_{\omega} \quad (28)$$

$$\lambda_y = \sum_{\omega} W_{\omega} (£_{\omega} - \mu_y) - (£_{\omega} - \mu_y)^T \quad (29)$$

## V. SIMULATION RESULTS

An interconnected microgrid test system made up of four interconnected MGs is considered as the case study to evaluate the performance of the proposed optimal scheduling model. The interconnected MG system has 68 sectionalizing switches along with 4 tie switches (one tie switch for each microgrid) as shown in Fig. 2. More detailed information about the test system's configuration is given in [21]. The interconnectivity within the system is illustrated by dotted lines. Table's I-III demonstrate the characteristics of distributed energy resources, energy storage, and adjustable loads in each MG.

Table IV shows the normalized hourly forecasted wind turbine generation, load demand, and electricity price during the 24-hour scheduling horizon where the peak load is 3,802.2 (kW). In order to simplify the problem, similar generation patterns with different capacities are considered for other three wind turbines. The optimization problem is solved using a 2.4-GHz personal computer with RAM 16-GB. In order to validate the effectiveness of the proposed problem, two cases are studied as follow:

***Case 1:*** Interconnected microgrid system optimal scheduling ignoring the microgrid distribution network constraints.

***Case 2:*** Interconnected microgrid system optimal scheduling considering the microgrid distribution network constraints.

Table I
CHARACTERISTICS OF DGs
(D: DISPATCHABLE, ND: NONDISPATCHABLE)

| | Unit | Type | Cost Coefficient ($/kWh) | Min-Max Capacity (kW) | Min Up/Down Time (h) | Ramp Up/ Down Rate (kW/h) |
|---|---|---|---|---|---|---|
| MG 1 | G1 | D | 0.154 | 800-3000 | 3 | 1500 |
| | G2 | D | 0.167 | 500-3500 | 3 | 1500 |
| | G3 | D | 0.182 | 800-2000 | 2 | 1000 |
| | G4 | ND | - | 0-1500 | - | - |
| MG 2 | G5 | D | 0.157 | 800-2000 | 3 | 1500 |
| | G6 | D | 0.134 | 800-2000 | 3 | 1500 |
| | G7 | D | 0.178 | 800-2000 | 3 | 1500 |
| MG 3 | G8 | D | 0.208 | 500-2500 | 3 | 1000 |
| | G9 | D | 0.195 | 700-3500 | 2 | 1500 |
| | G10 | ND | - | 0-2000 | - | - |
| MG 4 | G11 | D | 0.167 | 700-3000 | 3 | 1500 |
| | G12 | D | 0.184 | 500-2500 | 3 | 1000 |
| | G13 | ND | - | 0-1500 | - | - |

***Case 1:*** Table V shows results for the optimal DER scheduling during the 24-hours horizon. Based on the results, the most economic units are ON during the horizon to satisfy the load demand. For instance, in MG 1, only units 1 and 2 are ON which are the units with the lowest operation cost. In other words, the status of units is only based on the economic considerations. Table V presents adjustable loads scheduling for cases 1 and 2. Results indicate that the load scheduling is not affected by MG distribution network constraints.

TABLE II
CHARACTERISTICS OF THE ENERGY STORAGE IN INTERCONNECTED MICROGRID SYSTEM

| Storage | Capacity (kWh) | Min-Max Charging/Discharging Power (kW) | Min Charging/Discharging Time (h) |
|---|---|---|---|
| DES | 2000 | 50-200 | 5 |

TABLE III
CHARACTERISTICS OF ADJUSTABLE LOADS (S: SHIFTABLE, C: CURTAILABLE)

| Load | Type | Min-Max Capacity (kW) | Required Energy (kWh) | Initial Start/End Time (h) | Min Up Time (h) |
|---|---|---|---|---|---|
| L1(MG 1) | S | 0-120 | 480 | 11-14 | 1 |
| L2(MG 2) | S | 0-100 | 360 | 15-19 | 1 |
| L3(MG 3) | S | 20-80 | 240 | 16-19 | 1 |
| L4(MG 4) | C | 10-50 | 300 | 1-24 | 24 |
| L5(MG 5) | C | 30-90 | 400 | 13-24 | 12 |

TABLE IV
HOURLY FORECASTED WIND GENERATION (NORMALIZED), LOAD DEMAND (KW), AND ELECTRICITY PRICE ($/KWH)

| Time (h) | 1-6 | | | | | |
|---|---|---|---|---|---|---|
| Wind | 0.119 | 0.119 | 0.119 | 0.119 | 0.119 | 0.061 |
| Load | 0.800 | 0.805 | 0.810 | 0.818 | 0.830 | 0.910 |
| Price | 0.23 | 0.19 | 0.14 | 0.12 | 0.12 | 0.13 |
| Time (h) | 7-12 | | | | | |
| Wind | 0.119 | 0.087 | 0.119 | 0.206 | 0.385 | 0.394 |
| Load | 0.950 | 0.970 | 1.000 | 0.980 | 1.000 | 0.970 |
| Price | 0.13 | 0.14 | 0.17 | 0.22 | 0.22 | 0.22 |
| Time (h) | 13-18 | | | | | |
| Wind | 0.261 | 0.158 | 0.119 | 0.087 | 0.119 | 0.119 |
| Load | 0.950 | 0.900 | 0.905 | 0.910 | 0.930 | 0.900 |
| Price | 0.21 | 0.22 | 0.19 | 0.18 | 0.17 | 0.23 |
| Time (h) | 19-24 | | | | | |
| Wind | 0.087 | 0.119 | 0.087 | 0.087 | 0.061 | 0.041 |
| Load | 0.940 | 0.970 | 1.000 | 0.930 | 0.900 | 0.940 |
| Price | 0.21 | 0.22 | 0.18 | 0.17 | 0.13 | 0.12 |

***Case 2:*** In case 1, the interconnected configuration is not considered, thus the practical network limitations such as voltage limits and distribution line flow limits are not taken into account. In case 2, the proposed network models are added to the interconnected MG system optimal scheduling problem. Table VII represents the optimal DER scheduling during the operation horizon for the second scenario. According to table VII, the statuses of units are not only based on the economic considerations, but also they are affected by network limitations. For instance, in MG 1, all dispatchable units must commit to satisfy the load demand which shows the interconnectivity effect in the system.

Fig. 5 demonstrates the total operation cost for both case 1 and case 2 under deterministic and stochastic frameworks. According to this figure, by considering the stochastic framework, the total cost increases in both scenarios. This additional cost is paid to make a realistic analysis of the system with uncertainties.

TABLE VI

DER SCHEDULE IN CASE 1

| | Hours (1-24) | | | | | | | | | | | | | | | | | | | | | | | |
|---|---|---|---|---|---|---|---|---|---|---|---|---|---|---|---|---|---|---|---|---|---|---|---|---|
| G1 | 1 | 1 | 1 | 1 | 1 | 1 | 1 | 1 | 1 | 1 | 1 | 1 | 1 | 1 | 1 | 1 | 1 | 1 | 1 | 1 | 1 | 1 | 1 | 1 |
| G2 | 0 | 1 | 1 | 1 | 0 | 0 | 1 | 1 | 1 | 0 | 0 | 0 | 0 | 1 | 1 | 1 | 0 | 0 | 0 | 1 | 1 | 1 | 1 | 0 |
| G3 | 0 | 0 | 0 | 0 | 0 | 0 | 0 | 0 | 0 | 0 | 0 | 0 | 0 | 0 | 0 | 0 | 0 | 0 | 0 | 0 | 0 | 0 | 0 | 0 |
| G5 | 0 | 0 | 0 | 0 | 0 | 0 | 1 | 1 | 1 | 1 | 1 | 1 | 0 | 0 | 0 | 0 | 0 | 0 | 1 | 1 | 1 | 0 | 0 | 0 |
| G6 | 1 | 1 | 1 | 1 | 1 | 1 | 1 | 1 | 1 | 1 | 1 | 1 | 1 | 1 | 1 | 1 | 1 | 1 | 1 | 1 | 1 | 1 | 1 | 1 |
| G7 | 0 | 0 | 0 | 0 | 0 | 0 | 0 | 0 | 0 | 0 | 0 | 0 | 0 | 0 | 0 | 0 | 0 | 0 | 0 | 0 | 0 | 0 | 0 | 0 |
| G8 | 0 | 0 | 0 | 0 | 0 | 0 | 0 | 0 | 0 | 0 | 0 | 0 | 0 | 0 | 0 | 0 | 0 | 0 | 0 | 0 | 0 | 0 | 0 | 0 |
| G9 | 1 | 1 | 1 | 1 | 1 | 1 | 1 | 1 | 1 | 1 | 1 | 1 | 1 | 1 | 1 | 1 | 1 | 1 | 1 | 1 | 1 | 1 | 1 | 1 |
| G11 | 1 | 1 | 1 | 1 | 1 | 1 | 1 | 1 | 1 | 1 | 1 | 1 | 1 | 1 | 1 | 1 | 1 | 1 | 1 | 1 | 1 | 1 | 1 | 1 |
| G12 | 0 | 0 | 0 | 0 | 0 | 0 | 0 | 0 | 0 | 0 | 0 | 0 | 0 | 0 | 0 | 0 | 0 | 0 | 0 | 0 | 0 | 0 | 0 | 0 |
| DES | -1 | -1 | -1 | -1 | -1 | 1 | 1 | 1 | 1 | 1 | -1 | -1 | -1 | -1 | -1 | -1 | -1 | 1 | 1 | 1 | 1 | 1 | 1 | 1 |

TABLE V

ADJUSTABLE LOADS SCHEDULE IN CASES 1& 2

| | Hours (1-24) | | | | | | | | | | | | | | | | | | | | | | | |
|---|---|---|---|---|---|---|---|---|---|---|---|---|---|---|---|---|---|---|---|---|---|---|---|---|
| L1 | 0 | 0 | 0 | 0 | 0 | 0 | 0 | 0 | 0 | 0 | 1 | 1 | 1 | 1 | 0 | 0 | 0 | 0 | 0 | 0 | 0 | 0 | 0 | 0 |
| L2 | 0 | 0 | 0 | 0 | 0 | 0 | 0 | 0 | 0 | 0 | 0 | 0 | 0 | 0 | 1 | 1 | 1 | 1 | 1 | 0 | 0 | 0 | 0 | 0 |
| L3 | 0 | 0 | 0 | 0 | 0 | 0 | 0 | 0 | 0 | 0 | 0 | 0 | 0 | 0 | 0 | 1 | 1 | 1 | 1 | 0 | 0 | 0 | 0 | 0 |
| L4 | 1 | 1 | 1 | 1 | 1 | 1 | 1 | 1 | 1 | 1 | 1 | 1 | 1 | 1 | 1 | 1 | 1 | 1 | 1 | 1 | 1 | 1 | 1 | 1 |
| L5 | 0 | 0 | 0 | 0 | 0 | 0 | 0 | 0 | 0 | 0 | 0 | 0 | 1 | 1 | 1 | 1 | 1 | 1 | 1 | 1 | 1 | 1 | 1 | 1 |

TABLE VII

DER SCHEDULE IN CASE 2

| | Hours (1-24) | | | | | | | | | | | | | | | | | | | | | | | |
|---|---|---|---|---|---|---|---|---|---|---|---|---|---|---|---|---|---|---|---|---|---|---|---|---|
| G1 | 1 | 1 | 1 | 1 | 1 | 1 | 1 | 1 | 1 | 1 | 1 | 1 | 1 | 1 | 1 | 1 | 1 | 1 | 1 | 1 | 1 | 1 | 1 | 1 |
| G2 | 0 | 1 | 1 | 1 | 0 | 0 | 1 | 1 | 1 | 0 | 0 | 0 | 0 | 0 | 0 | 0 | 0 | 0 | 0 | 1 | 1 | 1 | 1 | 0 |
| G3 | 0 | 0 | 0 | 0 | 0 | 0 | 0 | 0 | 0 | 0 | 0 | 1 | 1 | 1 | 1 | 1 | 0 | 0 | 1 | 1 | 0 | 0 | 0 | 0 |
| G5 | 0 | 0 | 0 | 0 | 0 | 1 | 1 | 1 | 0 | 0 | 0 | 0 | 0 | 0 | 0 | 0 | 0 | 0 | 1 | 1 | 1 | 0 | 0 | 0 |
| G6 | 1 | 1 | 1 | 1 | 1 | 1 | 1 | 1 | 1 | 1 | 1 | 1 | 1 | 1 | 1 | 1 | 1 | 1 | 1 | 1 | 1 | 1 | 1 | 1 |
| G7 | 0 | 0 | 0 | 0 | 0 | 0 | 1 | 1 | 1 | 0 | 0 | 0 | 0 | 0 | 0 | 0 | 0 | 0 | 0 | 0 | 0 | 0 | 0 | 0 |
| G8 | 0 | 0 | 0 | 0 | 0 | 1 | 1 | 1 | 1 | 1 | 0 | 0 | 0 | 0 | 0 | 0 | 1 | 1 | 1 | 0 | 0 | 0 | 0 | 0 |
| G9 | 1 | 1 | 1 | 1 | 1 | 1 | 1 | 1 | 1 | 1 | 1 | 1 | 1 | 1 | 1 | 1 | 1 | 1 | 1 | 1 | 1 | 1 | 1 | 1 |
| G11 | 1 | 1 | 1 | 1 | 1 | 1 | 1 | 1 | 1 | 1 | 1 | 1 | 1 | 1 | 1 | 1 | 1 | 1 | 1 | 1 | 1 | 1 | 1 | 1 |
| G12 | 0 | 0 | 0 | 0 | 1 | 1 | 1 | 1 | 1 | 0 | 0 | 0 | 0 | 0 | 0 | 1 | 1 | 1 | 1 | 1 | 0 | 0 | 0 | 0 |
| DES | -1 | -1 | -1 | -1 | -1 | -1 | 1 | 1 | 1 | 1 | 1 | -1 | -1 | -1 | -1 | 1 | 1 | 1 | 1 | 1 | 0 | 0 | 0 | 0 |

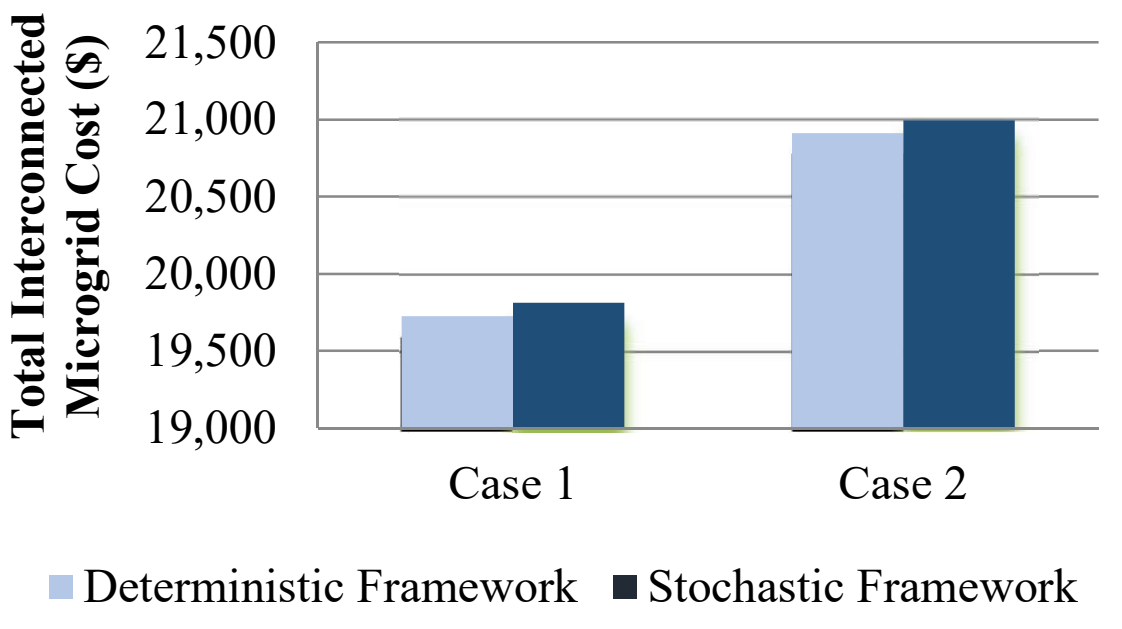

Fig. 5. Comparison of operation cost in different cases

## VI. CONCLUSION

This paper investigated the effects of interconnectivity within MGs from the energy management point of view. A newly developed stochastic framework based on the UT and a linear power flow is applied to the model solve the optimal scheduling of units in interconnected MGs under high uncertainties. The simulation results on an interconnected MG and through two different scenarios show that considering network constraints can affect the optimal scheduling of units. In other words, neglecting the physical limitations on power exchange between adjacent MGs can result in feeder congestion in the real analysis. From the uncertainty point of view, running the stochastic analysis has increased the interconnected MGs costs in both scenarios. This additional cost makes the analysis more reliable and practical compatible with the real-world limitations and forecast errors.

## REFRERENCES